\documentclass[a4paper,reqno]{amsart}

\usepackage{amsmath}
\usepackage{amsthm}
\usepackage{amsfonts}
\usepackage{amssymb}

\usepackage{color}

\usepackage{hyperref}

\hypersetup{
    colorlinks=false,
    pdfborder={0 0 0},
    pdfstartview={XYZ null null 1.00}
}

\numberwithin{equation}{section}

  \def\<{\langle}
  \def\>{\rangle}

\theoremstyle{plain}
  \newtheorem{theorem}{Theorem}[section]
  \newtheorem{proposition}[theorem]{Proposition}

\theoremstyle{definition}

  \newtheorem{remark}[theorem]{Remark}

\begin{document}

\title[Invariant sets and connecting orbits...]{Invariant sets and connecting orbits for nonlinear evolution equations at resonance}

\author{Piotr Kokocki}
\address{\noindent  Faculty of Mathematics and Computer Science \newline Nicolaus Copernicus University \newline  Chopina 12/18, 87-100 Toru\'n, Poland}
\email{pkokocki@mat.umk.pl}
\thanks{The researches supported by the NCN Grant no. 2011/01/N/ST1/05245}



\begin{abstract}
We study the problem of existence of orbits connecting stationary points for the
nonlinear heat and strongly damped wave equations being at resonance at infinity. The main
difficulty lies in the fact that the problems may have no solutions for general nonlinearity. To
address this question we introduce geometrical assumptions for the nonlinear term and use
them to prove index formulas expressing the Conley index of associated semiflows. We also
prove that the geometrical assumptions are generalizations of the well known Landesman-
Lazer and strong resonance conditions. Obtained index formulas are used to derive criteria 
determining the existence of orbits connecting stationary points.
\end{abstract}

\maketitle

\setcounter{tocdepth}{2}
%\tableofcontents

\section{Introduction}

Consider the following differential equations
\begin{align}\label{A-F-lam}
\dot u(t) & = - A u(t) + \lambda u(t) + F (u(t)),  \quad  t > 0 \\ \label{row-hyp1}
\ddot u(t) & = -A u(t) - c A \dot u(t) + \lambda u(t) + F(u(t)), \quad t > 0
\end{align}
where $A:D(A)\to X$ is a positive sectorial operator on a Banach space $X$ and $F:X^\alpha\to X$ is a continuous map on the fractional space $X^\alpha := D(A^\alpha)$, $\alpha\in (0,1)$.

Our objective is to study the existence of orbits connecting stationary points for these equations in the case of {\em resonance at infinity}, i.e. $$\mathrm{Ker}\, (\lambda I - A) \neq \{0\} \quad \text{and} \quad F \ \text{ is bounded}.$$
The main difficulty lies in the fact that, in the presence of resonance, the problem of existence of bounded orbits may not have solution for general nonlinearity $F$. This fact has been explained in detail in Remark \ref{rem-non-ex3}. Our aim is to overcome this difficulty by proving theorems determining the existence of orbits connecting stationary points for equations \eqref{A-F-lam} and \eqref{row-hyp1}, in the terms of appropriate geometrical assumptions imposed on the nonlinearity $F$. To this end we formulate below assumptions $(G1)$ and $(G2)$ and use them to prove the two main results: Theorems \ref{th-ind-orbi-hyp} and \ref{th-lan-laz}, which express the Conley index of the invariant set contained in sufficiently large ball in the terms of $(G1)$ and $(G2)$. These theorems are complement of results from \cite{MR798176, MR1992823}, where the parabolic equation with non-resonance condition at infinity is considered. 

Finally, we provide applications for particular partial differential equations. First of all, in Theorems \ref{lem-est2} and \ref{lem-est3}, we prove that if $F$ is a Nemytskii operator associated with a map $f:\Omega\times\mathbb{R}\to\mathbb{R}$, then the well-known Landesman-Lazer (see e.g. \cite{MR0267269}) and strong resonance conditions (see e.g. \cite{MR713209}) are actually particular cases of assumptions $(G1)$ and $(G2)$. 
Then we derive criteria determining the existence of orbits connecting stationary points for the heat and strongly damped wave equations.

\section{Spectral decomposition}

Let $A:D(A)\to X$ is a positive sectorial operator on a Banach space $X$ such that: \\[2pt]
\noindent\makebox[9mm][l]{$(A1)$}\parbox[t][][t]{120mm}{the operator $A$ has compact resolvents,}\\[2pt]
\noindent\makebox[9mm][l]{$(A2)$}\parbox[t][][t]{120mm}{there is a Hilbert space $H$ endowed with a scalar product $\langle\,\cdot\,, \,\cdot\,\rangle_H$ and norm $\|\cdot\|_H$ and a continuous injective map $i:X \hookrightarrow H$,}\\[2pt]
\noindent\makebox[9mm][l]{$(A3)$}\parbox[t][][t]{120mm}{there is a linear self-adjoint operator ${\widehat A}:H\supset D({\widehat A}) \to H$ such that $\mathrm{Gr}\,(A)\subset \mathrm{Gr}\,({\widehat A})$, where the graph inclusion is understood in the sense of product map $X \times X \xrightarrow{i\times i} H\times H$,}\\[5pt]

\begin{remark}\label{asd}
{\em One can prove (see e.g. \cite[Remark 3.1]{Kok3}) that the spectrum $\sigma(A)$ consists of the sequence of eigenvalues $\lambda_1 < \lambda_2 < \ldots < \lambda_i < \lambda_{i+1} < \ldots$ which is finite or $\lambda_i \to +\infty$ as $i\to +\infty$. Furthermore $\dim\mathrm{Ker}\,(\lambda_i I - A) < +\infty$ for $i\ge 1$. $\square$}
\end{remark}
In the following theorem we obtain spectral decomposition for the operator $A$.
\begin{theorem}{\em (\cite[Theorem 2.3]{Kok2})}\label{th:10}
If $\lambda = \lambda_k$ for some $k\ge 1$ is an eigenvalue of the operator $A$ and $X_0 := \mathrm{Ker}\, (\lambda I - A)$, then there are closed subspaces $X_+$, $X_-$ of $X$ such that $X = X_+\oplus X_-\oplus X_0$ and the following assertions hold: \\[2pt]
\noindent\makebox[5mm][r]{$(i)$} \parbox[t][][t]{120mm}{We have inclusions $X_-\subset D(A)$, $A(X_-)\subset X_-$, $A(X_+\cap D(A)) \subset X_+$ and furthermore $X_-$ is a finite dimensional space such that $X_-= \{0\}$ if $k=1$ and
$X_-=\bigoplus_{i=1}^{k-1} \mathrm{Ker}\,(\lambda_i I - A)$ if $k\ge 2$.}\\[4pt]
\noindent\makebox[5mm][r]{$(ii)$} \parbox[t][][t]{120mm}{If $A_+:X_+\supset D(A_+) \to X_+$ and $A_-:X_-\supset D(A_-) \to X_-$ are parts of the operator $A$ in $X_+$ and $X_-$, respectively, then $\sigma(A_+) = \{\lambda_i \ | \ i\ge k+1 \}$ and $\sigma(A_-) = \emptyset$ if $k=1$ and $\sigma(A_-) = \{\lambda_i \ | \ i=1,\ldots, k-1 \}$ if $k\ge 2$. }\\[4pt]
\noindent\makebox[5mm][r]{$(iii)$} \parbox[t][][t]{120mm}{The spaces $X_0$, $X_-$, $X_+$ are mutually orthogonal, that is, $\langle i (u_l),i(u_m)\rangle_H = 0$ for $l\neq m$, where $u_i\in X_i$ for $i\in\{0,-,+\}$.}
\end{theorem}

\section{Index formulas for invariant sets}

Let $A:D(A)\to X$ be a positive sectorial operator on a Banach space $X$ satisfying assumptions $(A1)$, $(A2)$, $(A3)$ and let $F:X^\alpha\to X$ be a continuous map on the fractional space $X^\alpha := D(A^\alpha)$, where $\alpha\in (0,1)$. Assume that \\[4pt]
\noindent\makebox[9mm][l]{$(F1)$}\parbox[t]{120mm}{for every $x\in X^\alpha$ there is an open neighborhood $V\subset X^\alpha$ of $x$ and constant $L > 0$ such that for $x_1,x_2\in V$ we have $\|F(x_1) - F(x_2)\|\le L \|x_1 - x_2\|_\alpha$,}\\[4pt]
\noindent\makebox[9mm][l]{$(F2)$}\parbox[t]{120mm}{there is a constant $m  > 0$ such that
$\|F(x)\| \le m$ for $x\in X^\alpha$,}\\[4pt]
\noindent\makebox[9mm][l]{$(F3)$}\parbox[t][][t]{120mm}{$F$ is {\em completely continuous}, that is, for any bounded set $V\subset X^\alpha$ the set $F(V)$ is relatively compact in $X$.}\\[3mm]

We recall that \emph{a mild solution} of the equation \eqref{A-F-lam} starting at $x$ is a continuous map $u\colon [0, +\infty) \to X^\alpha$ such that
\begin{equation*}
u(t) = e^{\lambda t}S_A(t)x + \int_0^t e^{\lambda (t - s)}S_A(t - s)F(u(s))\,d s \qquad\text{for}\quad t\ge 0.
\end{equation*}
It is the standard theory (see e.g. \cite[Theorem 3.3.3, Corollary 3.3.5]{MR610244}) that under the above assumptions, for any $x\in X^\alpha$, there is a unique mild solution $u(\,\cdot \, ; x)\colon [0,+\infty) \to X^\alpha$ of (\ref{A-F-lam}) such that $u(0; x) = x$. Therefore we are able to define \emph{the semiflow} $\Phi\colon [0,+\infty)\times X^\alpha \to X^\alpha$ for the equation (\ref{A-F-lam}) by $\Phi(t, x) := u(t; x)$ for $t\in[0,+\infty)$, $x\in X^\alpha$. Furthermore, note that equation \eqref{row-hyp1} may be written as
\begin{equation}\label{row-hyp}
    \dot w(t) = -{\bf A} w(t) + {\bf F}(w(t)), \qquad t > 0,
\end{equation}
where ${\bf A}:{\bf E}\supset D({\bf A})\to{\bf E}$ is a linear operator on ${\bf E}:=X^\alpha\times X$ given by
\begin{equation*}
\begin{aligned}
D({\bf A}) & :=\{(x,y)\in {\bf E} \ | \ x + c y\in D(A)\} \\
{\bf A} (x,y) & :=(-y,A(x + c y) - \lambda x)  \qquad \mathrm{for} \quad  (x,y)\in D({\bf A}),
\end{aligned}
\end{equation*}
and ${\bf F}:{\bf E}\to{\bf E}$ is given by ${\bf F}(x,y):=(0,F(x))$. Since $A$ is a sectorial operator we can prove that ${\bf A}$ is also a sectorial operator (see e.g. \cite{MR702424}, \cite{MR2571574}). Therefore, just as before we can define the semiflow ${\bf \Phi}\colon [0,+\infty)\times {\bf E} \to {\bf E}$ for the equation (\ref{row-hyp1}) by ${\bf \Phi}(t, x) := w(t; (x,y))$ for $t\in[0,+\infty)$, $(x,y)\in {\bf E}$, where $w(\,\cdot\,;(x,y)):[0,+\infty)\to {\bf E}$ is a solution for \eqref{row-hyp} starting at $(x,y)\in {\bf E}$.

We say that a map $u:\mathbb{R}\to X^\alpha$ (resp. $w:\mathbb{R}\to {\bf E}$) is {\em an orbit} provided $$ \Phi(t,u(s)) = u(t + s) \ \ (\text{resp. } {\bf \Phi}(t,w(s)) = w(t + s))\qquad\text{for}\quad t\ge 0, \ s\in\mathbb{R}.$$ We call the set $K\subset X^\alpha$ invariant, provided for every $x\in K$ there is an orbit $u$ for the semiflow $\Phi$ such that $u(0)\in K$ and $u(\mathbb{R})\subset K$. Similarly, $K\subset {\bf E}$ is an invariant if for every $(x,y)\in K$ there is a orbit $w$ for ${\bf \Phi}$ such that $w(0)\in K$ and $w(\mathbb{R})\subset K$.

Assume that $\lambda = \lambda_k$ for some $k\ge 1$ and consider a direct sum decomposition $X=X_0\oplus X_-\oplus X_+$ obtained in Theorem \ref{th:10}. Let $Q_1, Q_2, P:X\to X$ be the continuous projections onto $X_-$, $X_+$ and $X_0$, respectively. Define $$X_+^\alpha:= X^\alpha\cap X_+, \quad X_-^\alpha:= X^\alpha\cap X_- \quad\text{ and }\quad Q:= Q_1 + Q_2.$$

\begin{remark}\label{rem-non-ex3}
{\em If equation \eqref{A-F-lam} is at resonance at infinity then the problem of existence of compact orbits connecting stationary points may not have solution for general nonlinearity $F$. \\[4pt]
To see this it is enough to take $F(x) = y_0$ for $x\in X^\alpha$, where $y_0\in\mathrm{Ker}\, (\lambda I - A)\setminus\{0\}$. Indeed, if $u:\mathbb{R}\to X^\alpha$ is a bounded orbit, then $$u(t) = e^{\lambda (t - t')}S_A(t - t')u(t') + \int_{t'}^t e^{\lambda (t - \tau)}S_A(t - \tau) y_0 \, d\tau, \quad t > t'.$$ Since $\mathrm{Ker}\, (\lambda I - A) \subset \mathrm{Ker}\,(I - e^{\lambda t}S_A(t))$ for $t\ge 0$ it follows that $$u(t) = e^{\lambda (t - t')}S_A(t - t')u(t') + (t - t')y_0, \quad t > t',$$
and therefore, after acting by the operator $P$, we have $$Pu(t) = e^{\lambda (t - t')}S_A(t - t')Pu(t') + (t - t')Py_0 = Pu(t') + (t - t')y_0, \quad t > t',$$ and finally $Pu(h) = Pu(0) + h y_0$ for $h\ge 0$. This contradicts assumption that $u$ is bounded and proves our assertion. $\square$}
\end{remark}

To overcome this obstacles we introduce the following {\em geometrical assumptions:} \label{g1g2}
\begin{equation*}\leqno{(G1)}
\ \left\{\begin{aligned}
& \text{for every balls } B_1\subset X^\alpha_+\oplus X^\alpha_- \text{ and } B_2\subset X_0 \text{ there is } R > 0 \\ & \text{such that } \langle F(x + y), x\rangle_H > - \langle F(x + y), z\rangle_H \\
& \text{for } (y,z)\in B_1 \times B_2, \ x\in X_0 \text{ such that } \|x\|_H\ge R.
\end{aligned}\right.
\end{equation*}
and 
\begin{equation*}\leqno{(G2)}
\ \left\{\begin{aligned}
& \text{for every balls } B_1\subset X^\alpha_+\oplus X^\alpha_- \text{ and } B_2\subset X_0 \text{ there is } R > 0 \\ & \text{such that } \langle F(x + y), x\rangle_H < - \langle F(x + y), z\rangle_H \\
& \text{for } (y,z)\in B_1 \times B_2, \ x\in X_0 \text{ such that } \|x\|_H\ge R.
\end{aligned}\right.
\end{equation*}
With these assumptions we proceed to prove \emph{the index formulas for first and second order equations}. Assume that $\lambda = \lambda_k$ for $k\ge 1$, is an eigenvalue of $A$ and put $d_0:= 0$ and $d_l:=\sum_{i=1}^l \dim\mathrm{Ker}\,(\lambda_i I - A)$ for $l\ge 1$. The following index formula is a tool to determine the Conley index of the maximal invariant set contained in appropriately large ball for the equation \eqref{A-F-lam}.
\begin{theorem}{\em (\cite[Theorem 3.4]{Kok3})}\label{th-ind-orbi-hyp} 
There is a closed isolated neighborhood $N\subset X^\alpha$ such that, for $K:=\mathrm{Inv}\,(N,\Phi)$, the following statements hold: \\[3pt]
\noindent\makebox[4mm][r]{$(i)$} \parbox[t][][t]{120mm}{if condition $(G1)$ is satisfied, then $h(\Phi, K) = \Sigma^{d_k}$,}\\[3pt]
\noindent\makebox[4mm][r]{$(ii)$} \parbox[t][][t]{120mm}{if condition $(G2)$ is satisfied, then $h(\Phi, K) = \Sigma^{d_{k-1}}$.}
\end{theorem}
Here $h$ is the Conley index and $\Sigma^m$ is a homotopy type of pointed $m$-dimensional sphere (for more information see e.g. \cite{MR910097}). The following theorem is an analogous index formula for the equation \eqref{row-hyp1}.
\begin{theorem}{\em (\cite[Theorem 4.2]{Kok4})}\label{th-lan-laz} 
There is a closed isolated neighborhood $N\subset {\bf E}$, for $K:=\mathrm{Inv}\,(N,{\bf \Phi})$, we have the following assertions: \\[5pt] 
\noindent\makebox[4mm][r]{$(i)$} \parbox[t][][t]{120mm}{if condition $(G1)$ is satisfied, then $h({\bf \Phi}, K) = \Sigma^{d_k}$,}\\[3pt]
\noindent\makebox[4mm][r]{$(ii)$} \parbox[t][][t]{120mm}{if condition $(G2)$ is satisfied, then $h({\bf \Phi}, K) = \Sigma^{d_{k-1}}$.}
\end{theorem}

\section{Application to partial differential equations}

Assume that $\Omega\subset\mathbb{R}^n$ is an open bounded set with $C^\infty$ boundary and consider the following equations
\begin{align}\label{A-eps-res-a}
u_t(t, x) & = \Delta \, u(t,x) + \lambda  u(t,x) + f(x, u(t,x)) \\[5pt] \label{A-eps-res-ah}
u_{tt}(t, x) & = \Delta \, u(t,x) + c\Delta \, u_t(t,x) + \lambda  u(t,x) + f(x, u(t,x)).
\end{align}
where $c>0$ is a damping factor, $\lambda\in\mathbb{R}$ is a parameter, $\Delta$ is the Laplace operator with Dirichlet conditions, $f:\o\Omega\times\mathbb{R}\to\mathbb{R}$ is a continuous bounded map satisfying: \\[4pt]
\noindent\makebox[22pt][l]{$(E1)$} \parbox[t][][t]{120mm}{for every $R > 0$ there is $L > 0$ such that $|f(x,s_1) - f(x,s_2)| \le L(|s_1 - s_2|)$ for $x\in\Omega$ and $s_1,s_2\in\mathbb{R}$,}\\[4pt]
\noindent\makebox[22pt][l]{$(E2)$} \parbox[t][][t]{120mm}{$f$ is a map of class $C^1$ and there is a constant $\nu\in\mathbb{R}$ such that $\nu=D_s f(x,0)$ for $x\in\Omega$. Furthermore $f(x,0,0) = 0$ for $x\in\Omega$.} \\[3pt]

Put $\alpha\in (3/4, 1)$ and $p\ge 2n$ and let $X:=L^p(\Omega)$. Define $A_p: X\supset D(A_p)\to X$ as a linear operator given by $D(A_p) := W^{2,p}_0(\Omega)$ and $A_p \bar u := - \Delta \bar u$ for $\bar u\in D(A_p)$. \\[3pt]

\noindent\makebox[4mm][r]{$(1)$} \parbox[t][][t]{120mm}{It is known that $A_p$, $p\ge 2$, is a positive definite sectorial operator with compact resolvents and one can prove that $A_2$ is symmetric. Hence $(A1)$ is satisfied. }\\[3pt]
\noindent\makebox[4mm][r]{$(2)$} \parbox[t][][t]{120mm}{Take $H:=L^2(\Omega)$ with the standard inner product and norm. Since $\Omega$ is bounded and $p\ge 2$, we derive that $i:L^p(\Omega) \hookrightarrow L^2(\Omega)$ is a continuous embedding. In consequence we obtain assumption $(A2)$.}\\[3pt]
\noindent\makebox[4mm][r]{$(3)$} \parbox[t][][t]{120mm}{Using again the boundedness of $\Omega$, one can prove that for  $\widetilde{A}:= A_2$, the inclusion $A_p \subset \widetilde{A}$ is satisfied in the sense of the map $i\times i$. Therefore $(A3)$ holds.}\\[3pt]

By Remark \label{asd} it follows that the spectrum $\sigma(A)$ consists of the sequence of eigenvalues $\lambda_1 < \lambda_2 < \ldots < \lambda_i < \lambda_{i+1} < \ldots$ and $\dim\mathrm{Ker}\,(\lambda_i I - A) < +\infty$ for $i\ge 1$. The embedding theorem for  fractional spaces \cite[Theorem 1.6.1]{MR610244} implies that the inclusion $X^\alpha\subset C(\o\Omega)$ is continuous. Therefore we can define the {\em Nemytskii operator} $F\colon X^\alpha \to X$ given for every $\bar u\in X^\alpha$ by the formula $F(\bar u)(x) := f(x, \bar u(x))$ for $x\in\Omega$. \\[4pt]
Under the above assumptions the following assertions hold. \\[3pt]
\noindent\makebox[4mm][r]{$(i)$} \parbox[t][][t]{120mm}{One can prove that the map $F$ is continuous, bounded and satisfies assumptions $(F1)$, $(F2)$. Therefore, writing the equations \eqref{A-eps-res-a} and \eqref{A-eps-res-ah} in the abstract form \eqref{A-F-lam} and \eqref{row-hyp1}, respectively, we can associate with them semiflows $\Phi$ and ${\bf \Phi}$.}\\[4pt]
\noindent\makebox[4mm][r]{$(ii)$} \parbox[t][][t]{120mm}{One can also prove that $F$ is differentiable at $0$ and its derivative $DF(0)\in L(X^\alpha, X)$ is of the form $DF(0)[\bar u] = \nu \bar u$ for $\bar u\in X^\alpha$.}\\[4mm]
Now we proceed to examine what assumptions should satisfy the mapping $f$ in order to the associated Nemytskii operator $F$ meets the introduced earlier geometrical assumptions. We start with the following theorem which says that well-known \emph{Landesman-Lazer} conditions (see \cite{MR0267269}) are actually particular case of $(G1)$ and $(G2)$.
\begin{theorem}\label{lem-est2}
Assume that $f_+,f_-\colon \Omega \to \mathbb{R}$ are continuous functions such that
\begin{align*}
f_+(x) = \lim_{s \to +\infty} f(x,s) \quad\text{and}\quad f_-(x) = \lim_{s \to -\infty} f(x,s)  \qquad \mathrm{for} \quad  x\in\Omega.
\end{align*}
\makebox[4mm][r]{$(i)$} \parbox[t]{120mm}{If the condition 
$$\int_{\{\bar u>0\}} f_+(x) \bar u(x) \,d x  + \int_{\{\bar u<0\}} f_-(x) \bar u(x) \,d x > 0 \leqno{(LL1)}$$
is satisfied for $\bar u\in\mathrm{Ker}\,(\lambda I - A_p)\setminus\{0\}$, then $(G1)$ holds.}\\[7pt] 
\makebox[4mm][r]{$(ii)$} \parbox[t]{120mm}{If the condition
$$\int_{\{\bar u>0\}} f_+(x) \bar u(x) \,d x + \int_{\{\bar u<0\}} f_-(x) \bar u(x) \,d x < 0 \leqno{(LL2)}$$
is satisfied for $\bar u\in\mathrm{Ker}\,(\lambda I - A_p)\setminus\{0\}$, then $(G2)$ holds.} 
\end{theorem}
The following lemma proves that conditions $(G1)$ and $(G2)$ are also implicated by \emph{the strong resonance conditions}, studied for example in \cite{MR713209}.
\begin{theorem}\label{lem-est3}
Assume that there is a continuous function $f_\infty \colon \o\Omega \to \mathbb{R}$ such that
\begin{equation*}
f_\infty(x)  = \lim_{|s| \to +\infty} f(x,s)\cdot s \ \ \text{ for } \ x\in\o\Omega.
\end{equation*}
\makebox[4mm][r]{(i)} \parbox[t]{120mm}{Condition $(G1)$ is satisfied provided
\begin{equation*}\leqno{(SR1)}
\quad\left\{\begin{aligned}
& \text{there is } h\in L^1(\Omega) \text{ such that } f(x,s)\cdot s \ge h(x) \text{ for } (x,s)\in \Omega\times\mathbb{R} \\ &
\text{ and }  \int_\Omega f_\infty(x)\, d x > 0.
\end{aligned}\right.
\end{equation*}}\\
\makebox[4mm][r]{(ii)} \parbox[t]{120mm}{Condition $(G2)$ is satisfied provided
\begin{equation*}\leqno{(SR2)}
\quad\left\{\begin{aligned}
& \text{there is } h\in L^1(\Omega) \text{ such that } f(x,s)\cdot s \le h(x) \text{ for } (x,s)\in \Omega\times\mathbb{R} \\
& \text{ and } \int_\Omega f_\infty(x)\, d x < 0.
\end{aligned}\right.
\end{equation*}}
\end{theorem}  
The following theorem is a criterion on existence of orbits connecting stationary points with Landesman-Lazer type conditions.
\begin{theorem}\label{th-crit-ogrsrh}
Assume that $f_+,f_-\colon \Omega \to \mathbb{R}$ are continuous functions such that
\begin{equation*}
f_+(x) = \lim_{s \to +\infty} f(x,s) \quad\text{and}\quad  f_-(x) = \lim_{s \to -\infty} f(x,s) \ \ \text{for} \ x\in\Omega.
\end{equation*}
There is a nonzero compact orbit $u:\mathbb{R} \to X^\alpha$ of equation \eqref{A-eps-res-a} such that either $\lim_{t\to -\infty} u(t) = 0$ or $\lim_{t\to +\infty} u(t) = 0$, provided one of the following is satisfied: \\[2pt]
\makebox[4mm][r]{(i)} \parbox[t]{120mm}{$(LL1)$ holds and $\lambda_l < \lambda + \nu < \lambda_{l+1}$ where $\lambda_l \neq \lambda$;}\\[2pt]
\makebox[4mm][r]{(ii)} \parbox[t]{120mm}{$(LL1)$ holds and $\lambda + \nu < \lambda_1$;} \\[2pt]
\makebox[4mm][r]{(iii)} \parbox[t]{120mm}{$(LL2)$ holds and $\lambda_{l-1} < \lambda + \nu < \lambda_l$ where$\lambda \neq \lambda_l$, $l\ge 2$;}\\[2pt]
\makebox[4mm][r]{(iv)} \parbox[t]{120mm}{$(LL2)$ holds and $\lambda + \nu < \lambda_1$, where $\lambda \neq \lambda_1$.}
\end{theorem} 
In the proof of this theorem we need the following propositions.
\begin{proposition}{\em (see \cite[Theorem 3.5]{MR910097}, \cite[Proposition 4.3.3]{Kok5})}\label{con-ind-1}
If $\lambda + \mu \notin\sigma(A)$, then $h(\Phi,\{0\}) = h({\bf \Phi},\{0\}) = \Sigma^{b_l}$, where $b_l := 0$ if $\lambda + \mu < \lambda_1$ and $b_l := \sum_{i=1}^l \dim\mathrm{Ker}\,(\lambda_i I - A)$ if $\lambda_l < \lambda + \mu < \lambda_{l+1}$.
\end{proposition}
\begin{proposition}{\em (see \cite[Theorem 11.5]{MR910097})}\\\label{conn-orb}
\makebox[4mm][r]{(i)} \parbox[t]{120mm}{Let $K\subset X^\alpha$ be an isolated invariant set such that $0\in K$ and $\Phi(t,0) = 0$ for $t\ge 0$. If $h(\Phi,K) = \Sigma^m$, where $m\ge 0$ is an integer, $h(\Phi,\{0\}) \neq \o 0$ and $h(\Phi, \{0\}) \neq h(\Phi, K)$, then there is nonzero compact orbit $u:\mathbb{R} \to K$ of $\Phi$ such that either $\lim_{t\to -\infty}u(t) = 0$ or $\lim_{t\to +\infty}u(t) = 0$.}\\[4pt]
\makebox[4mm][r]{(ii)} \parbox[t]{120mm}{Let $K\subset {\bf E}$ be an isolated invariant set such that $0\in K$ and ${\bf \Phi}(t,0) = 0$ for $t\ge 0$. If $h({\bf \Phi},K) = \Sigma^m$, where $m\ge 0$ is an integer, $h({\bf \Phi},\{0\}) \neq \o 0$ and $h({\bf \Phi}, \{0\}) \neq h({\bf \Phi}, K)$, then there is nonzero compact orbit $w:\mathbb{R} \to K$ of ${\bf \Phi}$ such that either $\lim_{t\to -\infty}w(t) = 0$ or $\lim_{t\to +\infty}w(t) = 0$.}
\end{proposition} 
Now Theorem \ref{th-crit-ogrsrh} is a consequence of Theorems \ref{th-ind-orbi-hyp}, \ref{lem-est2} and Propositions \ref{con-ind-1}, \ref{conn-orb}. In a similar way, but using Theorem \ref{th-lan-laz} instead of Theorem \ref{th-ind-orbi-hyp} we can obtain the
following criterion for the existence of connecting orbits with strong resonance type conditions. 
\begin{theorem}\label{th-crit-ogrsrhaa}
Assume that there is a continuous function $f_\infty \colon \o\Omega \to \mathbb{R}$ such that
\begin{equation*}
f_\infty(x)  = \lim_{|s| \to +\infty} f(x,s)\cdot s \ \ \text{for} \ x\in\o\Omega.
\end{equation*}
Then there is a nonzero compact orbit $w:\mathbb{R} \to {\bf E}$ of equation \eqref{A-eps-res-ah} such that either $\lim_{t\to -\infty} w(t) = 0$ or $\lim_{t\to +\infty} w(t) = 0$, provided one of the following is satisfied: \\[2pt]
\makebox[4mm][r]{(i)} \parbox[t]{120mm}{condition $(SR1)$ holds and $\lambda_l < \lambda + \nu < \lambda_{l+1}$ where $\lambda_l \neq \lambda$;}\\[2pt]
\makebox[4mm][r]{(ii)} \parbox[t]{120mm}{condition $(SR1)$ holds and $\lambda + \nu < \lambda_1$;} \\[2pt]
\makebox[4mm][r]{(iii)} \parbox[t]{120mm}{condition $(SR2)$ holds and $\lambda_{l-1} < \lambda + \nu < \lambda_l$ where$\lambda \neq \lambda_l$, $l\ge 2$;}\\[2pt]
\makebox[4mm][r]{(iv)} \parbox[t]{120mm}{condition $(SR2)$ holds and $\lambda + \nu < \lambda_1$, where $\lambda \neq \lambda_1$.}
\end{theorem}

{\small

\parindent = 0 pt

\end{document}